\documentstyle [11pt] {article}

\input {amssym.def}
\input {amssym}
 
\title {Formal GNS Construction and WKB Expansion in 
        Deformation Quantization \\[3mm]
         {\normalsize Contribution to the Ascona Meeting on Deformation
          Theory, Symplectic Geometry, and Applications,
          \\ June 16-22, 1996}}

\author {{\bf M. Bordemann\thanks{mbor@phyq1.physik.uni-freiburg.de}~,
            \addtocounter{footnote}{2}
            S. Waldmann\thanks{waldman@phyq1.physik.uni-freiburg.de}}\\[3mm]
             Fakult\"at f\"ur Physik\\Universit\"at Freiburg\\
          Hermann-Herder-Str. 3\\79104 Freiburg i.~Br., F.~R.~G\\[3mm]
        }

\date{FR-THEP-96/21 \\[1mm]
      November 1996 \\[5mm]}

\newtheorem{theorem}{Theorem}

\newenvironment{prooof}{\begin{description}
                   \item[{\small {\bf Proof:}}] \small}{\hfill {\bf Q.E.D.}
                                                          \medskip
                                                       \end{description}}
 
\newtheorem{defi}{Definition}[theorem]
\newtheorem{prop}{Proposition}[theorem]
\newtheorem{lemma}[theorem]{Lemma}
\newtheorem{rem}{Remark}[theorem]
\newtheorem{cor}{Corollary}[theorem]

\newcommand{\bdef}{\begin{defi}}
\newcommand{\ede}{\end{defi}}
\newcommand{\bsat}{\begin{theorem}}
\newcommand{\esat}{\end{theorem}}
\newcommand{\bprop}{\begin{prop}}
\newcommand{\eprop}{\end{prop}}
\newcommand{\blem}{\begin{lemma}}
\newcommand{\elem}{\end{lemma}}
\newcommand{\brem}{\begin{rem}}
\newcommand{\erem}{\end{rem}}
\newcommand{\bcor}{\begin{cor}}
\newcommand{\ecor}{\end{cor}}
\newcommand{\bbew}{\begin{prooof}}
\newcommand{\ebew}{\end{prooof}}
\newcommand{\be}{\begin{equation}}
\newcommand{\ee}{\end{equation}}
\newcommand{\bea}{\begin{eqnarray}}
\newcommand{\eea}{\end{eqnarray}}
\newcommand{\beas}{\begin{eqnarray*}}
\newcommand{\eeas}{\end{eqnarray*}}
\newcommand{\ben}{\begin{enumerate}}
\newcommand{\een}{\end{enumerate}}
\newcommand{\lb}{\label}

\newcommand{\ra}{\rightarrow}

\newcommand{\f}{\frac}
\newcommand{\p}{\partial}

\newcommand{\Real} {\Bbb R}
\newcommand{\Complex}{\Bbb C}

\begin{document}

\maketitle

{\bf 1.} The concept of deformation quantization has been defined and 
exemplified in \cite{BFFLS78}. The existence of star-products on every 
symplectic manifolds has been established in \cite{DL83} and in \cite{Fed94}. 
This gives a reasonable physical picture of the noncommutative algebra
of quantum observables with built-in classical limit. However, the discussion
of a formal analogue of representations of the deformed algebra in some 
`Hilbert space' seems to have been restricted to examples in the literature
up to now. In \cite{BW96II} the authors have proposed how to construct
formal pre-Hilbert spaces for the deformed algebra in the same category
by means of a generalized Gel'fand-Naimark-Segal (GNS) construction.

We shall now briefly review this construction and shall apply this in
the next section to the Weyl star-product on the cotangent bundle of 
$\Real^n$ (which will give back the usual Schr\"odinger representation
together with the Weyl symmetrization rule; details thereof can be found
in \cite{BW96II}). The third section contains new material: 
We shall describe how to incorporate the usual WKB
expansion into the framework of star-products and GNS representations 
by means
of a certain positive linear functional on the deformed algebra having 
support on a projectable Lagrangean submanifold ${\rm graph}(dS)$ of $T^*\Real^n$.
The main trick is 
to use a suitable form of the star-exponential $e^{*\f{it}{\lambda}S}$.

Let $\Complex((\lambda))$ (resp. $\Real((\lambda))$) be the field of
complex (resp. real) Laurent series in $\lambda$ with finite principal
part (i.e. a finite number of negative powers in $\lambda$) 
and for any complex or real
vector space $V$ let $V((\lambda))$ denote the vector space of all formal
Laurent series (with finite principal part) in $\lambda$ with coefficients
in $V$. Let $(M,\omega)$ be a symplectic manifold equipped with a 
differential star-product $*$ satisfying $\overline{f*g}=\bar{g}*\bar{f}$.
Let $C^\infty(M)$ be the space of all smooth complex-valued functions on
$M$. Then $(C^\infty(M)((\lambda)),*)$ becomes an associative algebra over
the Laurent field $\Complex((\lambda))$. Since the field $\Real((\lambda))$
is an ordered field by using the real order of the lowest nonvanishing
coefficient, a $\Complex((\lambda))$-linear functional $\omega$
of $C^\infty(M)((\lambda))$ can be called {\em positive} iff 
$\omega(\bar{f}*f)$
is positive in the Laurent field $\Real((\lambda))$ for all 
$f\in C^\infty(M)((\lambda))$. As in the case of $C^*$-algebras it can
be shown that the so-called {\em Gel'fand ideal} 
$J_{\omega}:=\{f\in C^\infty(M)((\lambda))|\omega(\bar{f}*f)=0\}$ is a left
ideal of the algebra $C^\infty(M)((\lambda))$, that the quotient
${\cal H}_\omega:=C^\infty(M)((\lambda))/J_\omega$ (with canonical projection
denoted by $f\mapsto \Psi_f$) is a left module for 
$C^\infty(M)((\lambda))$, equipped with a positive definite 
$\Complex((\lambda))$-valued inner product defined by 
$\langle\Psi_f,\Psi_g\rangle:=\omega(\bar{f}*g)$, and that the 
GNS representation
$\pi_\omega(f)\Psi_g:=\Psi_{f*g}$ is a star-representation in the sense
that $\pi_\omega(f)$ and $\pi_\omega(\bar{f})$ are formal adjoints.
In \cite{BW96II} it was shown that on any K\"ahler manifold equipped
with the star-product of Wick type (see e.g. \cite {BW96I}) 
every delta-functional (evaluation
at a point) is a positve functional whose GNS-representation gives rise
to a formal Bargmann representation. 

\vspace{0.3cm}

{\bf 2.} Consider the cotangent bundle $\pi:T^*\Real^n\ra \Real^n$ of 
$\Real^n$
with its canonical symplectic structure, the standard Weyl product $*$,
and the standard volume form $d^nq:=dq^1\cdots dq^n$ on $\Real^n$. We define
\be
 C^\infty(T^*\Real^n)_{\Real^n}:=
    \{f\in C^\infty(T^*\Real^n)| {\rm ~supp}(f)\cap i(\Real^n)
                                   {\rm ~~is~compact.}  \}
\ee
where $i:\Real^n\ra T^*\Real^n$ is the canonical embedding of $\Real^n$
as the zero section. Then we define the $\Complex((\lambda))$-linear
functional $\omega_0:C^\infty(T^*\Real^n)\ra\Complex((\lambda))$ by
\be
      \omega_0(f):=\int_{\Real^n}\!d^nq~i^*f,
\ee
and let ${\cal S}:C^\infty(T^*\Real^n)((\lambda))\ra 
C^\infty(T^*\Real^n)((\lambda))$ be defined as 
${\cal S}:=\exp(-\f{i\lambda}{2}\Delta)$ where 
$\Delta:=\f{\p^2}{\p q^k\p p_k}$ then the following Lemma is found by 
integration by parts:
\blem
 For any two $f,g\in C^\infty(T^*\Real^n)_{\Real^n}((\lambda))$ we have
 \[
   \omega_0(f*g)=\int_{\Real^n}\!d^nq(i^*\bar{{\cal S}}f)(i^*{\cal S}g)
 \]
 whence $\omega_0$ is a positive linear functional. The Gel'fand ideal of
 $\omega_0$ is given by
 \[
   J_0:=\{ f\in C^\infty(T^*\Real^n)_{\Real^n}((\lambda))\; |\; 
   i^*{\cal S}f=0 \}.
 \]
\elem
In \cite{BW96II} we have shown that $J_0$ is the left ideal generated by
the $n$ momentum functions $p_1,\ldots,p_n$. Moreover we can determine the
GNS representation induced by $\omega_0$ explicitly using the map $\cal S$
and the fact that ${\cal S}\circ \pi^*=\pi^*$:
\bsat
 Let $f\in C^\infty(T^*\Real^n)((\lambda))$, 
 $g\in C^\infty(T^*\Real^n)_{\Real^n}((\lambda))$, and
 $\varphi,\psi\in C^\infty_0(\Real^n)((\lambda))$.
 \ben
  \item There is a canonical isomorphism 
   $\Phi:C^\infty(T^*\Real^n)_{\Real^n}((\lambda))/J_0\ra
    {\cal H}_0:=C^\infty_0(\Real^n)((\lambda))$
    given by $\Phi(\Psi_g):=i^*{\cal S}g$ whose inverse is simply given by 
    $\Phi^{-1}(\varphi)=\Psi_{\pi^*\varphi}$. The hermitean product on
    ${\cal H}_0$ induced by the pull-back of the GNS inner product with 
    respect to $\Phi$ is the standard $L^2$-inner product 
    $\langle\varphi,\psi\rangle=\int_{\Real^n}\!d^nq~\bar{\varphi}\psi$.
  \item $J_0$ is even a left ideal of $C^\infty(T^*\Real^n)((\lambda))$
    whence the GNS representation is defined for all its elements $f$
    and given by the following formula on elements of ${\cal H}_0$:
    \be \lb{5}
        \pi_0(f)\varphi=i^*{\cal S}(f*(\pi^*\varphi))
               =\sum_{r=0}^\infty\f{1}{r!}\left(\f{\lambda}{i}\right)^r
                 \f{\p^r({\cal S}f)}{\p p_{i_1}\cdots p_{i_r}}
                 \f{\p^r \varphi}{\p q^{i_1}\cdots q^{i_r}}.
    \ee
 \een
\esat
We had also shown in \cite{BW96II} that the representation (\ref{5})
exactly corresponds to the Weyl symmetrization rule of polynomials in
$q^1,\ldots,q^n,p_1,\ldots,p_n$.

\vspace {0.3cm}

{\bf 3.} 
We shall now use these results to formulate the WKB expansion for 
$T^*\Real^n$ in the 
framework of deformation quantization. Let $H\in C^\infty(T^*\Real^n)$
be a classical real-valued Hamiltonian and assume that there is a smooth
real-valued function $S:\Real^n\ra\Real$ satisfying the following (in general
nonlinear) first order partial differential equation, the so-called 
Hamilton-Jacobi equation, for some chosen real number $E$:
\be \lb{HamJac}
      H(dS(q))=E  {\rm ~~~for~all~}q\in\Real^n.
\ee
Let $L:={\rm graph}(dS)\subset T^*\Real^n$ be the graph or $dS$ which is known
to be a projectable Lagrangean submanifold of $T^*\Real^n$ (see e.g. 
\cite[p.28-31]{BaW95}). Equation 
(\ref{HamJac}) implies $(H-E)|_{L}=0$, and the local flow of the Hamiltonian
vector field of $H$ preserves $L$. 

We are now going to construct a suitable
GNS functional having support in $L$.
First, consider the
`Heisenberg equation' with respect to the Hamiltonian $\pi^*S$:
\be
      \f{df}{dt}(t)=\f{i}{\lambda}((\pi^*S)*f(t)-f(t)*(\pi^*S))
\ee
which can be shown to have a unique solution $f(t)=A_t(f_0)$ for each
initial function $f_0\in C^\infty(T^*\Real^n)((\lambda))$ where
$t\mapsto A_t$ is a one-parameter group of $\Complex((\lambda))$-linear
automorphisms of the $C^\infty(T^*\Real^n)((\lambda))$ commuting with
complex conjugation. It allows for the factorization
\be
       A_t=\Phi_t^*\circ T_t
\ee
where $\Phi_t$ is the Hamiltonian flow of $\pi^*S$ which consists of 
translations by $dS$ in the direction of the cotangent spaces, and 
$T_t=1+\sum_{r=1}^\infty \lambda^r T^{(r)}_t$
is a formal power series in $\lambda$ whose coefficients are differential
operators on $T^*\Real^n$.
Note that $\Phi_1$ is a diffeomorphism of the configuration space $\Real^n$
onto $L$.
Let $C^\infty(T^*\Real^n)_L$ be the space of all smooth complex-valued
functions on $T^*\Real^n$ whose support has compact intersection with $L$,
and let $\omega_1:C^\infty(T^*\Real^n)_L((\lambda))\ra\Complex((\lambda))$
be the following $\Complex((\lambda))$-linear functional (where $i_L$ denotes
the canonical injection $L\subset T^*\Real^n$):
\be
    \omega_1(f):=\int\!(\Phi_1^*d^nq)~i_L^*(T_{-1}f)~=~\omega_0\circ A_{-1}(f)
\ee
which obviously has support in $L$ and is easily seen to be positive 
since $A_{-1}$ is a real automorphism.
Moreover, the Gel'fand ideal of $\omega_1$ is simply given by $J_1=A_1J_0$,
and the GNS representation $\pi_1$ induced by $\omega_1$ can be constructed
as follows: since 
$C^\infty(T^*\Real^n)_L((\lambda))=
                   A_1(C^\infty(T^*\Real^n)_{\Real^n}((\lambda)))$
there is a well-defined unitary map $U:{\cal H}_0\ra{\cal H}_1$ between
the two formal GNS pre-Hilbert spaces given by 
$U(\Psi^{(0)}_f):=\Psi^{(1)}_{A_1f}$ whence the vector space ${\cal H}_1$
is isomorphic to $C^\infty_0(L)((\lambda))$. Observing that the Gel'fand ideal of
$\omega_1$ is again a left ideal in the algebra 
$C^\infty(T^*\Real^n)((\lambda))$ we get the following formula for the
GNS representation of $f\in C^\infty(T^*\Real^n)((\lambda))$:
\be
    \pi_1(f)=U\pi_0(A_{-1}f)U^{-1}.
\ee
We consider the time-independent Schr\"odinger equation
$\pi_1(H)\phi^{(1)}=E\phi^{(1)}$ for $\phi^{(1)}$ in the formal distribution
space ${\cal H}_1':=(C^\infty_0(L))'((\lambda))$. This is well-defined
since $\pi_1(H)$ is a series of differential operators. Via the map $U$,
this problem is related to the time-independent Schr\"odinger equation
\be \lb{SG}
    \pi_0(A_{-1}H)\phi=E\phi
\ee
where $\phi=\sum_{r=0}^\infty\lambda^r\phi_r$ and all the $\phi_r$ are
distributions in $C^\infty_0(\Real^n)'$. The above equation leads to a system
of partial differential equations for the $\phi_r$ given in the following
\bsat[Formal WKB Expansion] With the above notations we find that the 
       system of equations (\ref{SG}) is equaivalent to the following
       system of partial differential equations for $\phi$:
  \bea
    \lefteqn{(i^*\Phi_{-1}^*T_{-1}^{(1)}H)\phi_r
            -\f{i}{2}i^*\Delta((\Phi_{-1}^*H)(\pi^*\phi_r))
    +\f{i}{2}i^*\{ \Phi_{-1}^*H, \pi^*\phi_r \} =} \nonumber \\
   & &    -\sum_{\tiny \begin{array}{c}
                  a+b+c+d=r+1, \\
                  a,b,c,d>0,~r>d
              \end{array} }
              \f{(-i/2)^a}{a!}\f{(i/2)^b}{b!}
           \Delta^a(M_b(\Phi_{-1}^*T_{-1}^{(c)}H,\pi^*\phi_d)) \lb{15}
  \eea
   where $M_b$ are the bidifferential operators of the standard Weyl
   star product in $T^*\Real^n$, i.e. 
   $M_b(f,g):=
  (\p_{q^k}\p_{p'_k}-\p_{p_k}\p_{q'^k})^b(f(q,p)g(q',p'))|_{q=q',p=p'}$.
\esat

In the particular case of a Hamiltonian of the physical form
$H(q,p)=\sum_{i=1}^n(p_i)^2+V(q)$ for a smooth real-valued function $V$
on $\Real^n$ the system of equations (\ref{15}) reduces to the
well-known WKB transport equations (see e.g. \cite[13-14]{BaW95}) (which is
mainly due to the fact that $A_tH=\Phi_t^*H$ for $H$ at most quadratic in the
momenta) where $\phi_r:=0$ for negative $r$:
\be
    \f{\p^2S}{\p q^k\p q^k}\phi_r+2\f{\p S}{\p q^k}\f{\p \phi_r}{\p q^k}
                        =  i \f{\p^2\phi_{r-1}}{\p q^k\p q^k}.
\ee
Remarks:
\ben
 \item Since we assumed that the energy surface $H^{-1}(E)$ contains a graph
  of a one-form, namely $dS$, the physical situation in the above Theorem
  is in so far simplified as there are no classically forbidden regions or
  turning points. It is an interesting problem to perform the same GNS
  analysis for general cotangent bundles and more general Lagrangean 
  submanifolds incorporating Maslov corrections.
 \item The usual WKB phase function $e^{iS/\hbar}$ is recovered in the
  following way: assume first that $H$ is polynomial in the momenta.
  Then e. g. $\exp(*\f{it}{\lambda}(\pi^*S))*H$ is a well-defined formal 
  Laurent series in
  $1/\lambda$ converging for all $\lambda=\hbar\in\Real\setminus\{0\}$. 
  Hence
  $A_t(H)=\exp(*\f{it}{\lambda}(\pi^*S))*H*\exp(-*\f{it}{\lambda}(\pi^*S))=
  (\pi^*\exp(\f{it}{\lambda}S))*H*(\pi^*\exp(-\f{it}{\lambda}S))$ 
  converges for all 
  $\lambda=\hbar\in\Real\setminus\{0\}$. 
  Assuming that the power series $\phi$ converges
  for some $\hbar\in\Real\setminus\{0\}$ to a distribution, we readily 
  see that eqn (\ref{SG}) is equivalent to 
  $\pi_0(H)\exp(\f{it}{\hbar}S)\phi=\exp(\f{it}{\hbar}S)\phi$ since $\pi_0$
  is a representation.
\een

\section*{Acknowledgment}

\noindent We would like to thank M.~Flato, K.~Fredenhagen, J.~Huebschmann,
  H.~Reh\-ren, and A.~Wein\-stein for useful discussions.

\end{document}